\documentclass{article}
\usepackage{amssymb}
\usepackage{amsfonts}
\usepackage{amsmath}

\begin{document}

\title{On the Laws of Large Numbers in Possibility Theory}
\author{Sorin G. Gal\\
Department of Mathematics and Computer Science, \\
University of Oradea, \\
Universitatii 1, 410087, Oradea, Romania,\\
e-mail:\textit{galso@uoradea.ro}}
\date{}
\maketitle

\begin{abstract}
In this paper we obtain some possibilistic variants of the probabilistic laws of large numbers, different from those obtained by other authors, but very natural extensions of the corresponding ones in probability theory. Our results are based on the possibility measure and on the "maxitive" definitions for possibility expectation and possibility variance. Also, we show that in this frame, the weak form of the law of large numbers, implies the strong law of large numbers.
\end{abstract}



\section{Introduction}

It is well known the fact that possibility theory is an alternative theory to the probability theory, dealing with certain types of uncertainty and treatment of incomplete information (see, e.g., \cite{DubPrad} or \cite{Cooman}).
In the possibilistic models, all the probabilistic indicators (like expected value, variance, probability measure,
etc) are replaced with suitable possibilistic indicators.

Variants of the classical laws of large numbers in probability theory, were studied in the general setting of some non-additive set functions in a large number of papers, see, e.g., \cite{Fuller1}, \cite{Fuller2}, \cite{Hong}, \cite{Hong2}, \cite{Teran}, \cite{Teran2}, \cite{Cooman2}, \cite{Chen}.             .
They are, in essence, based on arithmetic means of the variables (fuzzy numbers), on some non-additive set valued functions including upper and lower probabilities, on very special concepts of possibilistic mean value, including those based on Choquet integrals.

In this paper, we consider a completely different approach. Thus, we use the "maxitive" concepts of possibilistic expectation and possibilistic variance as defined in \cite{DubPrad} (see also
\cite{Cooman}) and instead of the sum of a finite number of real functions, we consider the maximum of them.

It is worth noting that the idea of replacing for real functions the sum operator with the maximum operator, has been also proved very fruitful in constructing the so-called possibilistic (or max-product) approximation operators, see \cite{Gal-Pos}, \cite{CGOT}, \cite{BCG}. Among others, this approach allowed in the above mentioned works to introduce a possibilistic Feller's scheme, analogue
to the probabilistic Feller's scheme used in approximation by linear and positive operators.

Section 2 contains some preliminaries in the possibility theory, while in Section 3 we obtain some possibilistic variants of strong laws of large numbers.

\section{Preliminaries in Possibility Theory}

Firstly, in this section we present some known concepts and results in possibility theory, which will be useful in the next section. As it is easily seen, in fact they are the corresponding concepts for those in probability theory, like random variable, probability distribution, mean value, probability, so on.
For details, see, e.g., \cite{DubPrad} or \cite{Cooman}.

{\bf Definition 2.1.} Let $\Omega$ be a non-empty set.

(i) An application $X:\Omega\to \mathbb{R}$ will be called (random) variable.

(ii) A possibility distribution (on $\Omega$), is a function $\lambda : \Omega \to [0, 1]$, such that $\sup\{\lambda(s) ; s\in \Omega\}=1$.

(iii) The possibility expectation of a variable $X$ (on $\Omega$), with the possibility distribution $\lambda$, is defined by
$E_{sup}(X)=\sup_{s\in \Omega}X(s)\lambda(s)$. The possibility variance of $X$ is defined by $Var_{sup}(X)=\sup\{(X(s)-E_{sup}(X))^{2}\lambda(s) ; s\in \Omega\}$.

(iv) A possibility measure is a mapping $P:{\cal{P}}(\Omega)\to [0, 1]$, satisfying the axioms $P(\emptyset)=0$, $P(\Omega)=1$ and $P(\bigcup_{i\in I}A_{i})=\sup\{P(A_{i}) ; i\in I\}$ for all $A_{i}\in \Omega$, and any $I$, an arbitrary family of indices (if $\Omega$ is finite then obviously $I$ must be finite too). Note that if $A, B\subset \Omega$, satisfy $A\subset B$, then by the last property it easily follows that $P(A)\le P(B)$ and that $P(A\bigcup B)\le P(A)+P(B)$.

It is well-known (see, e.g., \cite{DubPrad}) that any possibility distribution $\lambda$ on $\Omega$, induces a possibility measure $P_{\lambda}:{\cal{P}}(\Omega)\to [0, 1]$, given by the formula $P_{\lambda}(A)=\sup \{\lambda(s) ; s\in A\}$, for all $A\subset \Omega$.

A key tool in our proofs is the following simple possibilistic analogue of the Chebyshev's inequality in probability theory.

{\bf Theorem 2.2.} (see, e.g., \cite{Gal-Pos}) {\it Let $\Omega$ be a non-empty set, $\lambda:\Omega \to [0, 1]$
and consider $X:\Omega \to \mathbb{R}$ be with the possibility distribution $\lambda$.
Then, for any  $r>0$, we have
$$P_{\lambda}(\{s\in \Omega ; |X(s)-E_{sup}(X)|\ge r\})\le \frac{Var_{sup}(X)}{r^{2}},$$
where $P_{\lambda}$ is the possibilistic measure induced by $\lambda$.}

This result was proved by Theorem 2.2 in \cite{Gal-Pos} for $\Omega$ discrete set, but analyzing its proof it is obvious that it remains valid for arbitrary non-empty sets too.

We also need the following concepts, known in fact for the more general case of non-additive set functions (see, e.g., \cite{Couso}, \cite{Denn}, \cite{WK1}).

{\bf Definition 2.3.} Let $X_{n}, X:\Omega\to \mathbb{R}$, $n\in \mathbb{N}$ and $P_{\lambda}$ be a possibility measure on $\Omega$.

(i) We say that $X_{n}$ converges to $X$ in the possibility measure $P_{\lambda}$, if for any $\varepsilon>0$, we have
$$\lim_{n\to \infty}P_{\lambda}(\{\omega\in \Omega ; |X_{n}(\omega)-X(\omega)|\ge \varepsilon \})=0.$$

(ii) We say that $X_{n}$ converges to $X$, almost everywhere with respect to the possibility measure $P_{\lambda}$, if
$$P_{\lambda}(\{\omega\in \Omega ; X_{n}(\omega)\not \to X(\omega)\})=0,$$
where the notation $X_{n}(\omega)\not \to X(\omega)$ includes the both cases when $(X_{n}(\omega))_{n}$ has not limit and when $(X_{n}(\omega))_{n}$ has a limit but it is different from $X(\omega)$.

{\bf Remark.} If $X_{n}$ converges to $X$, almost everywhere with respect to the possibility measure $P_{\lambda}$, then we have
$$P_{\lambda}(\{\omega\in \Omega ; \lim_{n\to \infty}X_{n}(\omega)=X(\omega)\})=1.$$
Indeed, denoting by $A=\{\omega\in \Omega ; X_{n}(\omega)\not \to X(\omega)\}$,
since any possibility measure is obviously subadditive we get
$$1=P_{\lambda}(\Omega)=P_{\lambda}(A\bigcup (\Omega\setminus A))\le P_{\lambda}(A)+P_{\lambda}(\Omega\setminus A)=
P_{\lambda}(\Omega\setminus A),$$
which obviously implies that $P_{\lambda}(\Omega\setminus A)=1$.

Also, we need the following auxiliary results.

{\bf Lemma 2.4.} {\it Let $P_{\lambda}$ be a possibility measure on $\Omega$ and $B_{n}\subset \Omega$, $n\in \mathbb{N}$. If $\lim_{n\to \infty}P_{\lambda}(B_{n})=0$, then
$\lim_{m\to \infty}\sup\{P_{\lambda}(B_{n}); n\ge m\}=0$, for all $n\in \mathbb{N}$.}

{\bf Proof.} Denoting $A_{m}=\sup\{P_{\lambda}(B_{n}); n\ge m\}$, $n\in \mathbb{N}$, it is immediate that the sequence of positive numbers $(A_{m})_{m\in \mathbb{N}}$ is non-increasing.

Let $\varepsilon >0$. By hypothesis, there exists $m_{0}\in \mathbb{N}$, such that $P_{\lambda}(B_{n})<\varepsilon$, for all $n\ge m_{0}$. This implies that $A_{m_{0}}=\sup\{P_{\lambda}(B_{n}) ; n\ge m_{0}\}<\varepsilon$ and since $(A_{m})_{m}$ is non-increasing, it follows $0\le A_{m}\le A_{m_{0}}<\varepsilon$, for all $m\ge m_{0}$.

In other words, $\lim_{m\to \infty}A_{m}=0$, which proves the lemma. $\hfill \square$

{\bf Lemma 2.5.} (Possibilistic Borel-Cantelli lemma) {\it Let $P_{\lambda}$ be a possibility measure on
$\Omega$. Denoting
$B_{\infty}=\bigcap_{m=1}^{\infty}\bigcup_{n=m}^{\infty}B_{n}$, if $\lim_{m\to \infty}\sup\{P_{\lambda}(B_{n}); n\ge m\}=0$, then $P_{\lambda}(B_{\infty})=0$.}

{\bf Proof.} Indeed, from the properties of the possibility measure $P_{\lambda}$ in Definition 2.1, (iv), it follows
$$P_{\lambda}(B_{\infty})\le P_{\lambda}(\bigcup_{n=m}^{\infty}B_{n})=\sup\{P_{\lambda}(B_{n}); n\ge m\}, \mbox{ for all } m\in \mathbb{N}.$$
Passing with $m\to \infty$, it follows $P_{\lambda}(B_{\infty})=0$.

It is useful to note here that in fact we can write
$$B_{\infty}=\{\omega\in \Omega ;  \mbox{ such that }  \omega \in B_{n} \mbox{ for infinitely many } n\}.$$
$\hfill \square$

{\bf Corollary 2.6.} {\it Let $X_{n}, X:\Omega\to \mathbb{R}$, $n\in \mathbb{N}$ and $P_{\lambda}$ be a possibility measure on $\Omega$.

If $X_{n}$ converges to $X$ in the possibility measure $P_{\lambda}$,
then $X_{n}$ converges to $X$, almost everywhere with respect to the possibility measure $P_{\lambda}$.}

{\bf Proof.} Let $\varepsilon>0$ be arbitrary. Denoting $B_{n}(\varepsilon)=\{\omega\in \Omega; |X_{n}(\omega)-X(\omega)|\ge \varepsilon\}$, since by hypothesis we have $\lim_{n\to \infty}P_{\lambda}(B_{n}(\varepsilon))=0$, Lemma 2.4 implies $\lim_{m\to \infty}\sup\{P_{\lambda}(B_{n}(\varepsilon)); n\ge m\}=0$, which by Lemma 2.5 implies $P_{\lambda}(B_{\infty}(\varepsilon))=0$, where
$$B_{\infty}(\varepsilon)=\{\omega\in \Omega ;  \mbox{ such that }  \omega \in B_{n}(\varepsilon) \mbox{ for infinitely many } n\}.$$
This means that for any  $\varepsilon>0$ and $\omega\in B_{\infty}(\varepsilon)$, we have that $X_{n}(\omega)\not \to X(\omega)$.

Now, let $(\varepsilon_{k})_{k}$ be a sequence with $\varepsilon_{k}\searrow 0$. It is clear that
$$\{\omega\in \Omega ; X_{n}(\omega)\not \to X(\omega)\}=\bigcup_{k=1}^{\infty}B_{\infty}(\varepsilon_{k}),$$
which implies
$$P_{\lambda}(\{\omega\in \Omega ; X_{n}(\omega)\not \to X(\omega)\})=\sup\{P_{\lambda}(B_{\infty}(\varepsilon_{k})) ;k\in \mathbb{N}\}=0.$$
Also, by the Remark after Definition 2.3, it follows that $P_{\lambda}(\{\omega\in \Omega ; X_{n}(\omega)\to X(\omega)\})=1$, which ends the proof of corollary. $\hfill \square$

{\bf Remarks.} 1) It is worth noting that different from what is happening in probability theory, in general, the almost everywhere convergence in Definition 2.3, (ii), does not imply the convergence in possibility measure.
A simple counterexample in this sense could be Example 7.6, p. 163 in \cite{WK1}. This is happening because of the following reasons. Indeed,
it is known that in probability theory, the validity of "a.e. convergence implies convergence in measure"  is based on the continuity from above of a probability measure. But it is easy to show that  while any possibility measure is continuous from below (that is
if $A_{1}\subset ... \subset A_{n}\subset A_{n+1}\subset ... \subset \Omega$ and  $A=\bigcup_{n=1}^{\infty}A_{n}$, then $P_{\lambda}(A)=\lim_{n\to \infty}P_{\lambda}(A_{n})$), in general it is not continuous from
above (that is by decreasing sequences of sets).

2) Corollary 2.6 will have as a consequence the fact that in this frame, the weak form of the law of large numbers, will always imply the strong law of large numbers.

\section{Laws of Large Numbers in Possibility Theory}

The following classical result in probability theory is well known.

{\bf Theorem 3.1.} (Kolmogorov, \cite{Kol}) {\it Let $(X_{k})_{k\in \mathbb{N}}$ be a sequence of independent random variables with finite variances $Var(X_{k})$ and denote $S_{n}=\sum_{k=1}^{n}X_{k}$. If $\sum_{k=1}^{\infty}\frac{Var(X_{k})}{k^{2}}<+\infty$, then $\frac{S_{n}-E(S_{n})}{n}\to 0$, almost sure in probability. Here $E(S_{n})$ denotes the probability expectation of $S_{n}$ and $Var(X_{k})$ denotes the probabilistic variance.}

In the same spirit of ideas, Petrov in \cite{Petrov1}, proved the following result.

{\bf Theorem 3.2.} {\it Let $(X_{k})_{k\in \mathbb{N}}$ be a sequence of independent random variables with finite variances $Var(X_{K})$ and denote $S_{n}=\sum_{k=1}^{n}X_{k}$. If $Var(X_{n})=O\left (\frac{n^{2}}{\Psi(n)}\right )$,
with a $\Psi$ positive, nondecreasing function satisfying $\sum_{n=1}^{\infty}\frac{1}{n \Psi(n)}<+\infty$,
then $\frac{S_{n}-E(S_{n})}{n}\to 0$, almost sure in probability.}

{\bf Remark.} Korchevsky proved in \cite{Kor1} that Theorem 3.2 is in fact a consequence of Theorem 3.1.

In this section, among others, we will obtain possibilistic variants of the above two results.

Thus, let $\Omega$ be a nonempty set and $P_{\lambda}:{\cal{P}}(\Omega)\to [0, 1]$ be a possibility measure generated by the possibility distribution $\lambda$.

Let $(X_{k})_{k\in \mathbb{N}}$ be a sequence of fuzzy variables. The main idea is to replace in the
above probabilistic laws of large numbers, the sum $\sum_{k=1}^{n}X_{k}(s)$ by the maximum $\max\{X_{1}(s), ..., X_{n}(s)\}, s\in \Omega$, the probabilistic variance and expectation with the their possibilistic variants in Definition 2.1, (iii) and the probability measure with the possibility measure in Definition 2.1, (iv).

The first possibilistic strong law of large numbers is the following.

{\bf Theorem 3.3.} {\it Let $\Psi(n)$, $n\in \mathbb{N}$, be such that $\Psi(n)>0$ for all $n\in \mathbb{N}$,
and $\Psi(n)\to \infty$. If the sequence $X_{k}:\Omega \to \mathbb{R}$, $k\in \mathbb{N}$, satisfies
$$\max\{Var_{sup}(X_{1}), ..., Var_{sup}(X_{n})\}\le C\cdot \left (\frac{n^{2}}{\Psi(n)}\right ), \mbox{ for all } n\in \mathbb{N},$$
with $C>0$ a constant independent of $n$, then denoting
$$M_{n}(\omega)=\max\{X_{1}(\omega), ..., X_{n}(\omega)\},$$
for $n\to \infty$ we have that $\frac{M_{n}-E_{sup}(M_{n})}{n}\to 0$, almost everywhere with respect to the measure of possibility $P_{\lambda}$. Also, it follows
$$P_{\lambda}\left (\left \{\omega\in \Omega ; \frac{M_{n}(\omega)-E_{sup}(M_{n})}{n}\not \to 0 \right \}\right )=0$$ and
$$P_{\lambda}\left (\left \{\omega\in \Omega ; \lim_{n\to \infty}\frac{M_{n}(\omega)-E_{sup}(M_{n})}{n}=0 \right \}\right )=1.$$}

{\bf Proof.} Firstly, we have
$$E_{sup}(M_{n})$$
$$=\sup\{M_{n}(\omega)\cdot \lambda(\omega); \omega\in \Omega\}=\sup\{\max\{X_{1}(\omega)\lambda(\omega), ..., X_{n}(\omega)\lambda(\omega)\}; \omega\in \Omega\}$$
$$=\max\{\sup\{X_{k}(\omega)\cdot \lambda(\omega) ; \omega\in \Omega\}; k=1, ..., n\}$$
$$=\max\{E_{sup}(X_{k}); k=1, ..., n\}.$$
Denoting $A_{n}(\omega)=\frac{M_{n}(\omega)}{n}$, the above calculation immediately implies
$$E_{sup}(A_{n})=\frac{\max\{E_{sup}(X_{k}); k=1, ..., n\}}{n}=\frac{E_{sup}(M_{n})}{n}.$$
In what follows we will use the inequality (see Lemma 11.3.1 in \cite{BCG}, p. 4.4.3)
$$|\max\{a_{i};i=1, ...,n\}-\max\{b_{i};i=1, ...,n\}|$$
$$\le \max\{|a_{i}-b_{i}|;i=1, ...,n\}, a_{i},b_{i}\in \mathbb{R}, i=1, ..., n.$$
Then, we get
$$Var_{sup}(A_{n})=$$
$$\sup\left \{\left (\frac{\max\{X_{1}(\omega), ...,X_{n}(\omega)\}}{n}-\frac{\max\{E_{sup}(X_{k}); k=1, ..., n\}}{n}  \right )^{2}\lambda(\omega); \omega\in \Omega \right \}$$
$$=\frac{1}{n^{2}}$$
$$\cdot \sup\left \{\left |\max\{X_{1}(\omega), ...,X_{n}(\omega)\}-\max\{E_{sup}(X_{k}); k=1, ..., n\}  \right |^{2}\lambda(\omega); \omega\in \Omega \right \}$$
$$\le \frac{1}{n^{2}}\sup\left \{\left (\max\{|X_{k}(\omega)-E_{sup}(X_{k})|; k=1, ..., n\}  \right )^{2}\lambda(\omega); \omega\in \Omega \right \}$$
$$=\frac{1}{n^{2}}\sup\left \{\max\{|X_{k}(\omega)-E_{sup}(X_{k})|^{2}; k=1, ..., n\}\lambda(\omega); \omega\in \Omega \right \}$$
$$=\frac{1}{n^{2}}\sup\left \{\max\{|X_{k}(\omega)-E_{sup}(X_{k})|^{2}\lambda(\omega); k=1, ..., n\}; \omega\in \Omega \right \}$$
$$=\frac{1}{n^{2}}\cdot \max\{Var_{sup}(X_{1}), ..., Var_{sup}(X_{n})\}\le \frac{1}{n^{2}}\cdot C \cdot \frac{n^{2}}{\Psi(n)}=\frac{C}{\Psi(n)}.$$
Applying now the Chebyshev's inequality in Theorem 2.2, for any $\varepsilon >0$ we immediately get
$$P_{\lambda}(\{\omega\in \Omega ; |A_{n}(\omega)-E_{sup}(A_{n})|\ge \varepsilon \})\le \frac{Var_{sup}(A_{n})}{\varepsilon^{2}}\le \frac{C}{\Psi(n) \varepsilon^{2}}.$$
This means that the sequence $Y_{n}(\omega)=\frac{M_{n}(\omega)-E_{sup}(M_{n})}{n}$, $n\in \mathbb{N}$, converges
to $0$ in the possibilistic measure $P_{\lambda}$.

Applying now Corollary 2.6 and the Remark after Definition 2.3, it immediately follows the conclusion of the theorem. $\hfill \square$

{\bf Remarks.}  1) If $\Psi(n)$, $n\in \mathbb{N}$, is such that $\Psi(n)>0$ for all $n\in \mathbb{N}$,
$\Psi(n)\nearrow \infty$, the sequence $\left (\frac{n^{2}}{\Psi(n)}\right )_{n}$ is non-decreasing in $n$ and
$(X_{k})_{k\in \mathbb{N}}$ satisfies $Var_{sup}(X_{k})\le C\cdot \left (\frac{k^{2}}{\Psi(k)}\right )$, for all $k\in \mathbb{N}$, then it is immediate that
$$\max\{Var_{sup}(X_{1}), ..., Var_{sup}(X_{n})\}\le C\cdot \left (\frac{n^{2}}{\Psi(n)}\right ), \mbox{ for all } n\in \mathbb{N},$$
which transform Theorem 3.3 into a kind of possibilistic analogue to Theorem 3.2.

2) Simple examples for the sequences $\Psi(n)$ in Theorem 3.3 and in the above Remark 1 are $\Psi(n)=n^{\delta}$ with $0<\delta \le 2$, or
$\Psi(n)=[\ln(n+1)]^{\delta}$ with $0<\delta \le 1$.

The second possibilistic strong law of large numbers is the following.

{\bf Theorem 3.4.} {\it Let $\delta\in (0, 2)$. If the sequence $X_{k}:\Omega \to \mathbb{R}$, $k\in \mathbb{N}$, satisfies
\begin{equation}\label{eq1}
\sup_{n\in \mathbb{N}}\frac{1}{n^{\delta}}\cdot \max\{Var_{sup}(X_{1}), ..., Var_{sup}(X_{n})\} = C < +\infty,
\end{equation}
then denoting
$$M_{n}(\omega)=\max\{X_{1}(\omega), ..., X_{n}(\omega)\},$$
we have that $Y_{n}(\omega):=\frac{M_{n}(\omega)-E_{sup}(M_{n})}{n}$ converges to $0$, almost everywhere with respect to $P_{\lambda}$. Also, it follows
$$P_{\lambda}\left (\left \{\omega\in \Omega ; \frac{M_{n}(\omega)-E_{sup}(M_{n})}{n}\not \to 0\right \}\right )=0,$$ and
$$P_{\lambda}\left (\left \{\omega\in \Omega ; \lim_{n\to \infty}\frac{M_{n}(\omega)-E_{sup}(M_{n})}{n}=0 \right \}\right )=1.$$}

{\bf Proof.} Denote $A_{n}=\frac{M_{n}}{n}$. By the proof of Theorem 3.3 and by hypothesis, we easily get
$$Var_{sup}(A_{n})\le \frac{\max\{Var_{sup}(X_{1}), ..., Var_{sup}(X_{n})\}}{n^{2}}$$
$$\le \frac{1}{n^{2-\delta}}\left (\frac{1}{n^{\delta}}\max\left \{Var_{\sup}(X_{1}), ..., Var_{\sup}(X_{n})\right \}\right )
\le \frac{1}{n^{2-\delta}}\cdot C, \mbox{ for all } n\in \mathbb{N}.$$
Applying now the Chebyshev's inequality in Theorem 2.2, for any $\varepsilon >0$ we immediately obtain
$$P_{\lambda}\left (\left \{\omega\in \Omega ; \left |\frac{M_{n}(\omega)-E_{sup}(M_{n})}{n}\right |\ge \varepsilon \right \}\right )\le \frac{Var_{sup}(A_{n})}{\varepsilon^{2}}\le \frac{C}{n^{2-\delta} \varepsilon^{2}}.$$
This means that $Y_{n}(\omega):=\frac{M_{n}(\omega)-E_{sup}(M_{n})}{n}$ converges to $0$ in the possibility measure $P_{\lambda}$, which by Corollary 2.6 and the Remark after Definition 2.3, it immediately implies the conclusion of the theorem. $\hfill \square$

{\bf Remark.} Theorem 3.4 could be considered as a possibilistic analogue of Theorem 2.6.1, p. 51 in the book \cite{Knill}, by replacing the condition (\ref{eq1}), by the stronger one
$$\sup_{n\in \mathbb{N}}\frac{1}{n^{\delta}}\cdot \sum_{k=1}^{n}Var_{sup}(X_{k})\le C < +\infty.$$

The third possibilistic strong law of large numbers can be considered as a kind of possibilistic analogue to the
Kolmogorov's result in Theorem 3.1.

{\bf Theorem 3.5.} {\it Let $\delta\in (0, 2)$. If the sequence $X_{k}:\Omega \to \mathbb{R}$, $k\in \mathbb{N}$, satisfies
\begin{equation}\label{eq2}
\sum_{k=1}^{\infty}\frac{Var_{sup}(X_{k})}{k^{\delta}}< +\infty,
\end{equation}
then denoting
$$M_{n}(\omega)=\max\{X_{1}(\omega), ..., X_{n}(\omega)\},$$
we have that $Y_{n}(\omega):=\frac{M_{n}(\omega)-E_{sup}(M_{n})}{n}$ converges to $0$, almost everywhere with respect to $P_{\lambda}$. Also, it follows
$$P_{\lambda}\left (\left \{\omega\in \Omega ; \frac{M_{n}(\omega)-E_{sup}(M_{n})}{n}\not \to 0\right \}\right )=0,$$ and
$$P_{\lambda}\left (\left \{\omega\in \Omega ; \lim_{n\to \infty}\frac{M_{n}(\omega)-E_{sup}(M_{n})}{n}=0 \right \}\right )=1.$$}

{\bf Proof.} Denoting $a_{k}=\frac{Var_{sup}(X_{k})}{k^{\delta}}$, by hypothesis (\ref{eq2}) it follows that $a_{k}\to 0$ as $k\to \infty$ and therefore there exists $C>0$ such that $0\le a_{k}\le C$, for all $k\in \mathbb{N}$. By the proof of Theorem 3.4 and keeping the notations there, we get
$$Var_{sup}(A_{n})\le \frac{\max\{Var_{sup}(X_{1}), ..., Var_{sup}(X_{n})\}}{n^{2}}$$
$$=\frac{1}{n^{2-\delta}}\cdot \max\left \{\frac{a_{k} k^{\delta}}{n^{\delta}} ; 1\le k\le n\right \}
\le \frac{1}{n^{2-\delta}}\cdot \max\left \{a_{k}; 1\le k\le n\right \}\le \frac{C}{n^{2-\delta}}.$$
From this point, the reasonings are identical with those in the proof of Theorem 3.4.
$\hfill \square$

{\bf Remark.} If to the hypothesis in Theorem 3.3, or in Theorem 3.4, or in Theorem 3.5 we add the hypothesis
that the sequence $(X_{k}(\omega))_{k}$ satisfies
$X_{k}(\omega)\ge 0$, for all $\omega\in \Omega$, $k\in \mathbb{N}$  and $E_{sup}(X_{k})=k\mu$, for all $k\in \mathbb{N}$, then we easily get that the sequence $\left (\frac{M_{n}(\omega)}{n}\right )_{n\in \mathbb{N}}$ converges as $n\to \infty$ to $\mu=E_{sup}(X_{1})$, almost everywhere in $\Omega$ with respect to the possibility measure $P_{\lambda}$.

{\bf Conclusion.} As we already have mentioned in Introduction, there are many papers devoted to the large laws of numbers for nonadditive set-functions. Their proofs are pretty technical due, in my opinion, to the use of arithmetical mean of variables. The natural replacement of the arithmetic mean of variables with their maximum (due to the "maxitive" property of possibility measures), makes the proofs of laws of large numbers so simple as in the probability theory. For this reason, we consider that the results in this paper can be considered as very natural extensions of those in probability theory o possibility theory.
The next interesting study would be to obtain central limit theorems in this frame.

\end{document}